# ASYMPTOTIC GLOBAL ROBUSTNESS IN BAYESIAN DECISION THEORY

By Christophe Abraham and Benoît Cadre

*ENSAM–INRA and Université Montpellier II*

In Bayesian decision theory, it is known that robustness with respect to the loss and the prior can be improved by adding new observations. In this article we study the rate of robustness improvement with respect to the number of observations $n$. Three usual measures of posterior global robustness are considered: the (range of the) Bayes actions set derived from a class of loss functions, the maximum regret of using a particular loss when the subjective loss belongs to a given class and the range of the posterior expected loss when the loss function ranges over a class. We show that the rate of convergence of the first measure of robustness is $\sqrt{n}$, while it is $n$ for the other measures under reasonable assumptions on the class of loss functions. We begin with the study of two particular cases to illustrate our results.

**1. Introduction.** In Bayesian analysis, choosing a prior distribution and choosing a loss function according to prior knowledge and preferences are difficult tasks. In practice, the decision maker usually chooses convenient approximations to the subjective prior and the subjective loss. The legitimacy of such approximations might be investigated by a sensitivity analysis of the results with respect to the approximations. This is the purpose of robust Bayesian analysis, which recently was overviewed by Ríos Insua and Ruggeri (2000). An interesting approach, called global robustness, proposes to replace a single prior distribution (resp. loss function) by a class of priors (resp. loss functions) and then to compute the range of the ensuing answers as the prior (resp. loss function) varies over the class.

Bayesians mainly focus on sensitivity to the prior distribution, although the final result can be drastically affected by the loss function. Moreover, Rubin (1987) showed that the loss function and the prior cannot be separated under a weak system of axioms for rational behavior. It is worth pointing out that robustness with respect to the prior can be expressed as a particular case of loss robustness. This is illustrated by the following example: the computation of the range of the posterior expectation when the prior density $p$ ranges over a class $\Gamma$ reduces to the computation of the range of the Bayes actions (i.e., decisions that minimize the posterior expected loss) when the loss function ranges over the class $\{l_2 p/p_0, p \in \Gamma\}$, where $l_2$ is the quadratic loss and $p_0$ is a fixed prior.

When robustness is lacking, Abraham (2001) showed it can be improved by adding new observations. It is of practical interest to know how many new observations are needed to achieve a given robustness. Herein we answer this question by investigating the asymptotic rate of convergence of three measures of posterior robustness. Because of the above remark, we focus on robustness with respect to the loss, since it provides a general framework including many prior robustness problems.

The asymptotic of global robustness measures (e.g., the range of posterior means or set probabilities) with respect to the prior has been investigated for particular classes (mainly $\varepsilon$-contamination classes) by Sivaganesan (1988), Pericchi and Walley (1991), Moreno and Pericchi (1993) and Ruggeri and Sivaganesan (2000). The local point of view has been studied by Gustafson and Wasserman (1995), Gustafson, Srinivasan and Wasserman (1996) and Sivaganesan (1996). For a recent account of the theory, refer to Sivaganesan (2000).

In Sections 4–6, we proceed with the study of three measures of posterior global robustness. Section 4 is devoted to the study of the Bayes actions set derived from a class of loss functions. We show that the Bayes actions set tends to a limit set with rate $\sqrt{n}$, where $n$ is the number of observations. In Section 5, we are concerned with the regret of choosing a decision associated with a particular loss function when the true loss function varies over a given class. We show that the rate of convergence of the supremum of the regrets is $\sqrt{n}$ or $n$, according to the class of loss functions. Section 6 deals with the range of the posterior expected loss, which has asymptotic rate $\sqrt{n}$ or $n$ as well. Section 2 provides two examples. For one of them, the above asymptotic rates are actually achieved for every finite $n$. In Section 3 we set up notation and terminology. In particular, we indicate that the posterior distribution can be calculated under misspecified models, that is, we contemplate that the observations are realizations from a convenient probability distribution with density $h_\sigma$ ($\sigma$ is the parameter), while the true distribution $Q$ may not correspond to $h_\sigma$ for all values of $\sigma$. Finally, we compile some auxiliary results in Section 8.

**2. Examples.** In this section, we present two examples based on tractable classes of loss functions. Such classes have already been considered in Martín, Ríos Insua and Ruggeri (1998) and Abraham and Daurés ((1999, 2000)).

2.1. *Squared-error loss.* Whereas squared-error loss is frequently used to approximate nearly symmetric loss functions [Berger (1985)], it is of practical interest to investigate robustness with respect to variations around this loss. It is also of theoretical interest because it makes the calculations relatively simple.

The set $\Theta$ of parameters and the set $\mathcal{D}$ of decisions are both assumed to be $\mathbb{R}$. Fix $0 < k_1 < k_2$, depending on the incomplete information on the true loss, and define $U \colon \Theta \times \mathcal{D} \to \mathbb{R}^+$ as

$$(2.1) \qquad U(\sigma, d) = (k_2\{d \geq \sigma\} + k_1\{d < \sigma\})l_0(\sigma, d),$$

where $\{C\}$ denotes the usual indicator function of $C$ and $l_0(\sigma, d) = 0.5(d - \sigma)^2$ denotes the convenient loss chosen by the decision maker. Define $L$ by interchanging $k_1$ and $k_2$ in the definition of $U$. Let $D_{01}l$ stand for the derivative of $l \colon \Theta \times \mathcal{D} \to \mathbb{R}^+$ with respect to $d$ and introduce the class $\mathcal{F}$ of loss functions $l \colon \Theta \times \mathcal{D} \to \mathbb{R}^+$ such that for all $\sigma \in \Theta$, $D_{01}l(\sigma, \cdot)$ is continuous, $l(\sigma, \sigma) = 0$ and $D_{01}L \leq D_{01}l \leq D_{01}U$.

Assume that $X_1, \ldots, X_n$ are independent and identically distributed from a normal $N(\mu, \lambda^{-1})$ distribution, where the variance $\lambda^{-1}$ is known. Take a $N(\mu_0, \lambda_0^{-1})$ prior. The posterior $\pi_n$ is then normal $N(\mu_n, \lambda_n^{-1})$ with $\mu_n = (\lambda_0\mu_0 + \lambda(X_1 + \cdots + X_n))/\lambda_n$ and precision $\lambda_n = \lambda_0 + n\lambda$. Denoting, for all $l \in \mathcal{F}$, $d_l^n$ as a minimizer of $l^n(\cdot) = \int_\Theta l(\sigma, \cdot)\pi_n(d\sigma)$, elementary calculations show that $U^n$ and $L^n$ admit only one minimizer given by

$$d_U^n = \mu_n + r_1/\sqrt{\lambda_n} \quad \text{and} \quad d_L^n = \mu_n + r_2/\sqrt{\lambda_n},$$

where $r_2 < 0 < r_1$ are constants depending on $k_1$ and $k_2$.

Let us now investigate the computation of the three measures of posterior robustness. Since, by Abraham and Daurés (1999), $\{d_l^n, l \in \mathcal{F}\} = [d_U^n, d_L^n]$, the diameter of $\{d_l^n, l \in \mathcal{F}\}$ is equal to $(r_1 - r_2)/\sqrt{\lambda_n}$, which gives the first measure of robustness. Write now $\operatorname{reg}_l^n(d) = l^n(d) - \inf_\mathcal{D} l^n$ for the posterior regret. By the definition of $\mathcal{F}$, if $d_2 \geq d_1$, we have for all $l \in \mathcal{F}$,

$$l^n(d_2) - l^n(d_1) = \int_{d_1}^{d_2} D_{01}l^n(t)\,dt = \int_{d_1}^{d_2}\int_\Theta D_{01}l(\sigma, t)\pi_n(d\sigma)\,dt.$$

Hence, we deduce that

$$(2.2) \qquad \sup_{l \in \mathcal{F}} \operatorname{reg}_l^n(d) = \max\{\operatorname{reg}_U^n(d), \operatorname{reg}_L^n(d)\}.$$

Let $d_0^n = \mu_n$ be the Bayes rule associated with the squared-error loss function $l_0$. After some calculations, we obtain that for some constants $c_1$ and $c_2$,

$$U^n(d_0^n) - U^n(d_U^n) = c_1/\lambda_n \quad \text{and} \quad L^n(d_0^n) - L^n(d_L^n) = c_2/\lambda_n,$$

and hence

$$\sup_{l \in \mathcal{F}} \operatorname{reg}_l^n(d_0^n) = \max(c_1, c_2)/\lambda_n,$$

which gives the second measure of robustness. Finally, if $S = k_2 l_0$ and $I = k_1 l_0$, then $I, S \in \mathcal{F}$ and $I \leq l \leq S$ for all $l \in \mathcal{F}$. Then if we write $\operatorname{ran}^n(d) = \sup_{l \in \mathcal{F}} l^n(d) - \inf_{l \in \mathcal{F}} l^n(d)$ for the range of the posterior expected loss, we obviously have

$$\operatorname{ran}^n(d_0^n) = S^n(d_0^n) - I^n(d_0^n)$$
$$= 0.5(k_2 - k_1)/\lambda_n,$$

hence the third measure of robustness.

We emphasize that the constants $r_1$, $r_2$, $c_1$ and $c_2$ can be numerically computed and that similar calculations can be done with different functions $U$, $L$ and $l_0$. As a conclusion, we proved that, for the class $\mathcal{F}$, the speed of convergence of the diameter of $\{d_l^n, l \in \mathcal{F}\}$ is $\sqrt{n}$, while the speed of convergence of the posterior regret and the range of the posterior expected loss are $n$.

2.2. *The dam construction problem.* Following Ulmo and Bernier (1973), the economical consequence of constructing a dam $d$ meters high is the sum of the cost construction and the cost due to a potential flood, $10d + 100(H - d)\{H > d\}$, where $H$ is the peak water level. Note that the consequence is a random variable. Assuming that $H$ is exponentially distributed with density $h_\sigma(x) = \sigma e^{-\sigma x}$ and taking the expectation yields the loss

$$l_0(\sigma, d) = 10d + 100\sigma^{-1} \exp(-d\sigma).$$

A similarly constructed utility function can be found in Berger [(1985), page 58]. The loss $l_0$ can be viewed as a convenient approximation to the true loss. Let us proceed similarly to Section 2.1 to study the robustness of the Bayes action. Consider the class $\mathcal{F}$ of functions $l$ such that $D_{01}L \leq D_{01}l \leq D_{01}U$. Whereas the minimum of $l_0(\sigma, \cdot)$ is obtained when $d\sigma = \log 10$, we define

$$U(\sigma, d) = (\Phi(d\sigma - \log 10) + 0.5)\ l_0(\sigma, d)$$

and

$$L(\sigma, d) = (1.5 - \Phi(d\sigma - \log 10))\ l_0(\sigma, d),$$

where $\Phi$ denotes the cumulative distribution function of $N(0, 1)$. Let $d_l^n$ and $d_0^n$ be generic notation for the Bayes actions associated with the loss functions $l$ and $l_0$, respectively. It can be proved that $U(\sigma, \cdot)$ and $L(\sigma, \cdot)$ are convex functions with a unique minimizer. Thus, the set of Bayes actions is still $[d_U^n, d_L^n]$ and the largest posterior regret can be calculated by (2.2). The posterior distribution is derived from $n$ independent observations with density $h_\sigma$ and a reference prior $\pi(\sigma) = \sigma^{-1}$ $(\pi_n \sim \text{Gamma}(n, \sum_{i=1}^n X_i))$. We simulated $n = 100$ observations with respect to $h_{0.5}$ and computed numerically $\sum_{i=1}^{100} x_i = 193.6$, $d_U^n = 2.7$, $d_L^n = 7.7$, $d_0^n = 4.5$ and $\sup_{l \in \mathcal{F}} \text{reg}_l^n(d_0^n) = 19.5$. Thus, the optimal dam size is somewhere between 2.7 and 7.7 m, and using the optimal decision associated with $l_0$ gives an excess posterior loss less than 19.5. Can we get more precise results by adding new observations? Sections 4 and 5 answer in the negative. Indeed, Theorem 4.1 applied to $\mathcal{L} = \{U, L\}$ shows that the range of the optimal sizes approaches $d_L^\theta - d_U^\theta$ with rate $\sqrt{n}$, where $\theta$ is the true value of the parameter $\sigma$, and $d_U^\theta$ and $d_L^\theta$ are the minimizers of $U(\theta, \cdot)$ and $L(\theta, \cdot)$. From the data we can guess $\theta$ to be about 0.5 (because $1/\bar{x} = 0.51$) and deduce that $d_L^\theta - d_U^\theta$ is around 5 by numerical computation of $d_L^\theta$ and $d_U^\theta$ for $\theta = 0.5$. Since $d_L^n - d_U^n = 5$, we cannot expect to improve the result. Note that the class $\mathcal{F}$ is large since, even when $\theta$ is given, it is only known that the optimal size is somewhere between $d_U^\theta$ and $d_L^\theta$. Also note that if we had chosen a class $\mathcal{F}$ such that $d_U^\theta = d_L^\theta$, the range of the optimal sizes could have been arbitrarily reduced by adding observations [see Abraham (2001) for a description of the limit of the Bayes actions set]. Similarly, we know from Theorem 5.1 that the largest posterior regret approaches $\max\{\text{reg}_L^\theta(d_0^\theta), \text{reg}_U^\theta(d_0^\theta)\}$, which remains about 20, where $d_0^\theta$ denotes the minimizer of $l_0(\theta, \cdot)$.

## 3. Preliminaries and notation.

3.1. *The model.* Let $X = (X_1, X_2, \ldots)$ be a sample sequence of independent and identically distributed random variables defined on some measurable space $(\mathcal{X}_0, \mathcal{B}_0)$, where $\mathcal{B}_0$ denotes the Borel $\sigma$-field of $\mathcal{X}_0$. In the sequel $Q$ refers to the joint distribution on $(\mathcal{X}, \mathcal{B})$ of the sequence $X$, where $\mathcal{X} = \mathcal{X}_0^{\mathbb{N}}$ and $\mathcal{B}$ denotes the Borel $\sigma$-field of $\mathcal{X}$.

We introduce the family of probability densities $\{h_\sigma, \sigma \in \Theta\}$ with respect to some $\sigma$-finite measure $\mu$ on $(\mathcal{X}_0, \mathcal{B}_0)$, where the parameter space $\Theta$ is $\mathbb{R}^k$ with Borel $\sigma$-field $\mathcal{B}_\Theta$. Note that the model may be misspecified since we do not assume that $Q$ corresponds to any of the densities $h_\sigma$. For technical reasons, we make the additional assumption that $(\sigma, x_0) \to h_\sigma(x_0)$ is $\mathcal{B}_\Theta \otimes \mathcal{B}_0$ measurable.

From now on, we fix a prior distribution $\pi$ on $(\Theta, \mathcal{B}_\Theta)$. The existence of the posterior distribution for misspecified models was studied by Berk (1970). For simplicity, we assume that the posterior distribution $\pi_n$, defined for all $A \in \mathcal{B}_\Theta$ by

$$\pi_n(A) = \int_A \prod_{i=1}^n h_\sigma(X_i)\pi(d\sigma) \bigg/ \int_\Theta \prod_{i=1}^n h_\sigma(X_i)\pi(d\sigma),$$

does exist $Q$-almost surely.

We assume the model $h_\sigma$ to be regular enough so that the maximum likelihood estimate $\theta_n$ is asymptotically normal [i.e., for some $\theta \in \Theta$, $\sqrt{n}(\theta_n - \theta)$ converges in distribution to a normal random variable $Z_\theta$] and the posterior distribution concentrates around the true value of the parameter as $n \to \infty$. The precise assumptions M on the model are given in the beginning of Section 8. Sufficient conditions for the existence and the asymptotic normality of $\theta_n$ (i.e., assumption M1) with misspecified models were given by White (1982) for the case when $\Theta$ is compact. Moreover, Abraham and Cadre (2002) studied the concentration of $\pi_n$ around the true value of the parameter; see also Strasser (1976) when the model is correctly specified. More precisely, both works give sufficient conditions so that M2–M4 hold.

3.2. *The basic class of loss functions.* For simplicity, let $\mathcal{D} = \mathbb{R}^p$ be the decision space. In the sequel a loss function is defined to be a function $l \cdot \Theta \times \mathcal{D} \to \mathbb{R}^+$ such that $l(\cdot, d)$ is measurable for each $d \in \mathcal{D}$ and $l(\sigma, \cdot)$ is twice continuously differentiable for each $\sigma \in \Theta$.

If $\sigma_i$ (resp. $d_i$) denotes the $i$th component of $\sigma \in \Theta$ (resp. $d \in \mathcal{D}$), we write, when they exist,

$$D_{01}l = \left(\frac{\partial l}{\partial d_i}\right)_{i=1,\ldots,p}, \qquad D_{10}l = \left(\frac{\partial l}{\partial \theta_i}\right)_{i=1,\ldots,k},$$

$$D_{02}l = \left(\frac{\partial^2 l}{\partial d_i \partial d_j}\right)_{i,j=1,\ldots,p}, \qquad D_{20}l = \left(\frac{\partial^2 l}{\partial \theta_i \partial \theta_j}\right)_{i,j=1,\ldots,k},$$

$$D_{11}l = \left(\frac{\partial^2 l}{\partial d_i \partial \theta_j}\right)_{\substack{i=1,\ldots,p \\ j=1,\ldots,k}},$$

where $i$ and $j$ stand for the row index and the column index, respectively.

In this article a class $\mathcal{L}$ of loss functions is said to be locally $\pi$-dominated if, for all $d \in \mathcal{D}$, there exist a function $g \in L_1(\pi)$ which is bounded on a neighborhood of $\theta$, and an open ball $B(d, r)$ with center $d$ and radius $r > 0$ such that

$$\sup_{l \in \mathcal{L}} \sup_{t \in B(d,r)} \|D_{0\gamma}l(\sigma, t)\| \le g(\sigma), \qquad \sigma \in \Theta, \ \gamma = 0, 1, 2,$$

with the notation $D_{00}l = l$. Here and in the sequel $\|a\|$ denotes the maximum of the absolute values of the coordinates of a vector or a matrix $a$ with real coefficients. Thus, a locally $\pi$-dominated class is also locally $\pi_n$-dominated on the event $\{\int g(\sigma)\pi_n(d\sigma) < \infty\}$, the probability of which tends to 1 when $n \to \infty$ by Lemma 8.1. Since this article deals with convergence in probability and in distribution, we may restrict our attention to the elements of this set.

To shorten notation, we write $l^n(d) = \int_\Theta l(\sigma, d)\pi_n(d\sigma)$ as the expectation of $l(\cdot, d)$ with respect to $\pi_n$. Note that in a locally $\pi$-dominated class differentiation and integration can be inverted, and we let

$$D_{0\gamma} l^n(d) = \int_\Theta D_{0\gamma} l(\sigma, d)\pi_n(d\sigma), \qquad \gamma = 1, 2.$$

Furthermore, if $D_{0\gamma} l(\sigma, \cdot)$ is continuous for $\pi$-almost all $\sigma$, $D_{0\gamma} l^n$ is continuous as well.

3.3. *The Bayes action process.* Since, for each loss function $l$, $l^n(d)$ is a measurable function of $x$ and a continuous function of $d$, it is possible, for each $x \in \mathcal{X}$ such that $\arg\min_{d \in \mathcal{D}} l^n(d) \neq \varnothing$, to select a minimizing decision $d_l^n(x)$ in such a manner that the function $x \mapsto d_l^n(x)$ is $\mathcal{B}$ measurable [Rockafellar and West (1998), Theorem 14.37]. The decision $d_l^n$ is called the Bayes action associated with the loss $l$.

We use the outer probability theory to avoid strong assumptions on $\mathcal{L}$ that ensure the measurability of $(d_l^n)_{l \in \mathcal{L}}$. We denote by $Q^*$ the outer probability derived from $Q$, by $Y_n \xrightarrow{Q^*} Y$ the convergence in outer probability and by $Y_n \rightsquigarrow Y$ the weak convergence (with respect to $Q^*$) of $Y_n$ to $Y$. For more details about outer probability, refer to van der Vaart and Wellner (1996).

Throughout this article $\mathcal{L}$ denotes a locally $\pi$-dominated class of loss functions such that the outer probability that $\arg\min_{d \in \mathcal{D}} l^n(d) = \varnothing$ for some $l \in \mathcal{L}$ is zero. We then define a Bayes actions process to be a family $(d_l^n)_{l \in \mathcal{L}}$ of minimizing decisions. We equip the space of functions from $\mathcal{L}$ into the space of matrices with real coefficients with the supremum norm.

**4. Asymptotic of the Bayes actions process.** This section is devoted to the study of the Bayes actions process. To get asymptotic results, it is necessary to put some restrictions on $\mathcal{L}$. We assume throughout that $\mathcal{L}$ satisfies the following properties [recall that $\theta$ is fixed (see Section 3.1)]:

1a. For every $l \in \mathcal{L}$, $\arg\min l(\theta, \cdot) = \{d_l^\theta\}$.
1b. There exists a neighborhood $V_\theta$ of $\theta$ such that, for all $l \in \mathcal{L}$, $D_{01} l(\cdot, d_l^\theta)$ is continuously differentiable on $V_\theta$.
1c. $\sup_{l \in \mathcal{L}} \|D_{11} l(\theta, d_l^\theta)\| < \infty$, $\sup_{l \in \mathcal{L}} \|D_{02} l(\theta, d_l^\theta)\| < \infty$ and $\inf_{l \in \mathcal{L}} |\det D_{02} l(\theta, d_l^\theta)| > 0$.
1d. The families $\{D_{11} l(\cdot, d_l^\theta)|_{V_\theta}, l \in \mathcal{L}\}$, $\{D_{02} l(\cdot, d_l^\theta)|_{V_\theta}, l \in \mathcal{L}\}$ and $\{l(\cdot, d)|_{V_\theta}, l \in \mathcal{L}, d \in K\}$ are equicontinuous at $\theta$ for any compact $K \subset \mathcal{D}$.

Let $B(c, r)$ be generic notation for an open ball with center $c$ and radius $r > 0$.

1e. For every $\eta > 0$, there exists $\rho_\eta \in L_1(\pi)$ with $\sup_{\sigma \in V_\theta} \rho_\eta(\sigma) \to_{\eta \to 0} 0$ and such that for all $\sigma \in \Theta$ we have

$$\sup_{l \in \mathcal{L}} \sup_{d \in B(d_l^\theta, \eta)} \|D_{02} l(\sigma, d) - D_{02} l(\sigma, d_l^\theta)\| \leq \rho_\eta(\sigma).$$

1f. There exist $r > 0$ and a compact set $K \subset \mathcal{D}$ such that

$$\sup_{l \in \mathcal{L}} \inf_{d \in K} l(\theta, d) < \inf_{l \in \mathcal{L}} \inf_{\sigma \in B(\theta, r)} \inf_{d \in K^c} l(\sigma, d).$$

1g. For every $\eta > 0$,

$$\kappa(\eta) = \inf_{l \in \mathcal{L}} \inf_{d \in B^c(d_l^\theta, \eta)} [l(\theta, d) - l(\theta, d_l^\theta)] > 0.$$

The homogeneity of $\mathcal{L}$ is ensured by conditions 1b–1e. From 1f we prove that the Bayes actions remain in a compact set (Lemma 8.2). Let us illustrate the assumptions by the following examples.

EXAMPLE 4.1 (Prior robustness). Let $\Gamma$ be a class of densities with respect to (w.r.t.) the Lebesgue measure $m$ on $\mathbb{R}$ and assume $\pi$ has a positive density $w_0$ w.r.t. $m$. Consider the class $\mathcal{L}$ of functions $l(\sigma, d) = (d - a(\sigma))^2 w(\sigma)/w_0(\sigma)$ with $w \in \Gamma$. For instance, we take $a(\sigma) = \sigma$ or $a(\sigma) = \{\sigma \in S\}$ whether we are interested in the posterior expectation or the posterior probability of a set $S$. For simplicity, let us choose $a(\sigma) = \sigma$. Assume that $w_0$ and each $w \in \Gamma$ are continuously differentiable on a neighborhood $V_\theta$ of $\theta$. If furthermore $\sup_{w \in \Gamma} \sup_{\sigma \in V_\theta} w(\sigma) < \infty$, $\sup_{w \in \Gamma} \sup_{\sigma \in V_\theta} |w'(\sigma)| < \infty$ and $\inf_{w \in \Gamma} \inf_{\sigma \in V_\theta} w(\sigma) > 0$, assumptions 1a–1g are fulfilled.

Classes as in Example 4.1 include density band classes, mixture classes and $\varepsilon$-contamination classes with adequate conditions. [Conditions on the $\varepsilon$-contamination class are those used by Sivaganesan (1996).]

EXAMPLE 4.2. Consider the case $\Theta = \mathcal{D} = \mathbb{R}$. Assume that $\int_\Theta |\sigma|^p \times \pi(d\sigma) < \infty$ and let $g : \mathbb{R} \to [0, \infty)$ be a polynomial of degree $p$. Consider the class $\mathcal{G}$ of three times differentiable non-negative functions $f$ such that $|f^{(3)}(t)| \leq g(t)$. Assume further that $f$ is decreasing on $(-\infty, 0]$ and increasing on $[0, \infty)$ with a unique minimizer at 0 and that there exists $M > 0$ such that $\sup_{f \in \mathcal{F}} f(0) < \infty$, $\sup_{f \in \mathcal{F}} f(0) < \inf_{f \in \mathcal{F}} \inf_{|t| > M} f(t)$ and $0 < \inf_{f \in \mathcal{F}} f''(0) \leq \sup_{f \in \mathcal{F}} f''(0) < \infty$. Then the class $\mathcal{L}$ of loss functions $l(\sigma, d) = f(d - \sigma)$, $f \in \mathcal{G}$, satisfies every assumption of Section 3.2 and 1a–1g of Section 4.

This example includes, for instance, parametric classes (with Linex losses) and $\varepsilon$-contamination classes with adequate conditions [for definitions and examples of classes of loss functions, refer to Ríos Insua and Ruggeri (2000)].

To shorten notation, we write $\varphi(l)$ instead of $[D_{02} l(\theta, d_l^\theta)]^{-1} D_{11} l(\theta, d_l^\theta)$.

THEOREM 4.1. Under the assumptions M:

(i) $\sqrt{n} \sup_{l \in \mathcal{L}} \| (d_l^n - d_l^\theta) + \varphi(l)(\theta_n - \theta) \| \xrightarrow{Q_*^n} 0.$

(ii)  $\sqrt{n}(d_l^n - d_l^\theta)_{l \in \mathcal{L}} \rightsquigarrow (\varphi(l) Z_\theta)_{l \in \mathcal{L}}$.

(iii)  $\sqrt{n} \sup_{l \in \mathcal{L}} \|d_l^n - d_l^\theta\| \rightsquigarrow \sup_{l \in \mathcal{L}} \|\varphi(l) Z_\theta\|$.

From a robust point of view it is of interest to know the rate of convergence of the Bayes actions set with respect to the Hausdorff metric $h$. Let $\mathcal{A} = \{d_l^\theta, l \in \mathcal{L}\}$ and $\mathcal{A}^n = \{d_l^n, l \in \mathcal{L}\}$. Recall that $h(\mathcal{A}^n, \mathcal{A}) < \delta$ if and only if every point in $\mathcal{A}$ is within distance $\delta$ of at least one point in $\mathcal{A}^n$ and vice versa. Thus, $h(\mathcal{A}, \mathcal{A}^n) \le \sup_{l \in \mathcal{L}} \|d_l^n - d_l^\theta\|$ and, by Theorem 4.1,

$$\sqrt{n}/u_n h(\mathcal{A}, \mathcal{A}^n) \overset{Q^*}{\to} 0$$

for any sequence of positive numbers such that $u_n \to \infty$, thus improving the main result in Abraham (2001). Clearly, the same result holds if $h(\mathcal{A}, \mathcal{A}^n)$ is replaced by (diameter $\mathcal{A}^n$ − diameter $\mathcal{A}$). Assuming moreover that $\mathcal{D} = \Theta = \mathbb{R}$ and $d_l^\theta = d^\theta$ is independent of $l \in \mathcal{L}$, we get from Theorem 4.1,

$$\sqrt{n} \text{ diameter } \mathcal{A}^n \rightsquigarrow \sup_{l \in \mathcal{L}} (\varphi(l) Z_\theta) - \inf_{l \in \mathcal{L}} (\varphi(l) Z_\theta).$$

EXAMPLE 4.1 (continued).  Assume that for some $\overline{w} \in \Gamma$ with $\int_\Theta \sigma^2 \times \overline{w}(\sigma) \, d\sigma < \infty$ we have $w \le \overline{w}$ for all $w \in \Gamma$. The class $\mathcal{L}$ is then $\pi$-dominated. Write $\tilde{l}^n(d) = \int_\Theta (d - \sigma)^2 w_n(d\sigma)$, where $w_n$ is the posterior distribution derived from the prior density $w$, and denote by $\tilde{\mathcal{A}}^n$ the set of posterior expectations. Since $\varphi(l) = -1$ and $\tilde{\mathcal{A}}^n = \mathcal{A}^n$, we deduce from above that $\sqrt{n}$ diameter $\tilde{\mathcal{A}}^n \rightsquigarrow 0$.

EXAMPLE 4.2 (continued).  Since $\varphi(l) = -1$, we have $\sqrt{n}$ diameter $\mathcal{A}^n \rightsquigarrow 0$.

PROOF OF THEOREM 4.1.  Recall that integration and differentiation can be interchanged in a locally $\pi$-dominated class. By definition of $d_l^n$, $0 = D_{01} l^n(d_l^n)$. For $s \in [0, 1]$ write $t_{l,s}^n = d_l^\theta + s(d_l^n - d_l^\theta)$. By Taylor's formula we have

$$0 = \sqrt{n} D_{01} l^n(d_l^\theta) + \sqrt{n} \int_0^1 D_{02} l^n(t_{l,s}^n)^t (d_l^n - d_l^\theta) \, ds$$

$$= \sqrt{n}(D_{01} l^n(d_l^\theta) - D_{11} l(\theta, d_l^\theta)(\theta_n - \theta))$$

$$\quad + \left[ \int_0^1 D_{02} l^n(t_{l,s}^n)^t \, ds \right] \sqrt{n}(d_l^n - d_l^\theta) + D_{11} l(\theta, d_l^\theta) \sqrt{n}(\theta_n - \theta)$$

$$= \alpha_n(l) + A_n(l) \sqrt{n}(d_l^n - d_l^\theta) - R_n(l)$$

with evident definitions of $\alpha_n(l)$, $A_n(l)$ and $R_n(l)$. By Theorem 8.1 the supremum when $l$ ranges over $\mathcal{L}$ of $\alpha_n(l)$ tends to 0 in outer probability. Then (i) is straightforward from Lemmas 8.4 and 8.6. By Slutsky's lemma and M1, (i) gives (ii). Taking into account the continuity of the application $z \to \sup_{l \in \mathcal{L}} \|z(l)\|$, where $z$ is a function from $\mathcal{L}$ to $\mathbb{R}^k$, we easily deduce (iii) from (ii). □

**5. Posterior regret.** Let $l_0 \in \mathcal{L}$. From now on we think of $l_0$ as a convenient approximation of the true loss. For simplicity of notation we write $d_0^\theta$ and $d_0^n$ instead of $d_{l_0}^\theta$ and $d_{l_0}^n$. We let $\mathcal{S}_0 \subset \mathcal{L}$ be a class which satisfies the following conditions (recall that $V_\theta$ and $\rho_\eta$ were defined by 1b and 1e):

2a. For every $l \in \mathcal{S}_0$, $l(\cdot, d_0^\theta)$ is continuously differentiable on $V_\theta$.

2b. For every $\eta > 0$ and $\sigma \in \Theta$, we have

$$\sup_{l \in \mathcal{S}_0} \sup_{d \in B(d_0^\theta, \eta)} \|D_{01} l(\sigma, d) - D_{01} l(\sigma, d_0^\theta)\| \leq \rho_\eta(\sigma).$$

2c. The families $\{D_{01} l(\cdot, d_0^\theta)|_{V_\theta}, \, l \in \mathcal{S}_0\}$ and $\{D_{10} l(\cdot, d_0^\theta)|_{V_\theta}, \, l \in \mathcal{S}_0\}$ are equicontinuous at $\theta$.

2d. $\sup_{l \in \mathcal{S}_0} \|D_{01} l(\theta, d_0^\theta)\| < \infty$ and $\sup_{l \in \mathcal{S}_0} \|D_{10} l(\theta, d_0^\theta)\| < \infty$.

Similarly, the class $\mathcal{S} \subset \mathcal{L}$ is defined by replacing $d_0^\theta$ by $d_l^\theta$ and $\mathcal{S}_0$ by $\mathcal{S}$ in conditions 2a–2d. In the remainder of this section we restrict our attention to a class of loss functions $\mathcal{L}_1 \subset \mathcal{S} \cap \mathcal{S}_0$.

For every $l \in \mathcal{L}_1$ and every $d \in \mathcal{D}$, write

$$\mathrm{reg}_l^n(d) = l^n(d) - \inf_{d \in \mathcal{D}} l^n(d) \quad \text{and} \quad \mathrm{reg}_l^\theta(d) = l(\theta, d) - \inf_{d \in \mathcal{D}} l(\theta, d).$$

This section is devoted to the study of the posterior regret process for the decision $d_0^n$ associated with the convenient loss $l_0$. This measure of robustness was used by Berger (1984).

THEOREM 5.1. *Under the assumptions* M,

$$\sqrt{n}(\mathrm{reg}_l^n(d_0^n) - \mathrm{reg}_l^\theta(d_0^\theta))_{l \in \mathcal{L}_1}$$

$$\rightsquigarrow ([-D_{01} l(\theta, d_0^\theta)^t \varphi(l) + D_{10} l(\theta, d_0^\theta)^t - D_{10} l(\theta, d_l^\theta)^t] Z_\theta)_{l \in \mathcal{L}_1}.$$

Taking into account the continuity of the application $z \to \sup_{l \in \mathcal{L}_1} \|z(l)\|$ defined on the functions from $\mathcal{L}_1$ to $\mathbb{R}^k$, we deduce from Theorem 5.1 the asymptotic bound for every $u \in \mathbb{R}$,

$$\limsup_n Q^* \left( \sqrt{n} \sup_{l \in \mathcal{L}_1} |\mathrm{reg}_l^n(d_0^n) - \mathrm{reg}_l^\theta(d_0^\theta)| \geq u \right) \leq Q \left( \sup_{l \in \mathcal{L}_1} |M_l| \geq u \right),$$

where $(M_l)_{l \in \mathcal{L}_1}$ is the limit process that appears in Theorem 5.1. The above inequality provides information on the value of $n$ that we need to obtain an arbitrarily robust analysis. For instance, choose $\alpha$ arbitrarily small and $u \in \mathbb{R}$ so that the right-hand term is less than $\alpha$. Then with probability greater than $1 - \alpha$, the posterior regret $\mathrm{reg}_l^n(d_0^n)$ associated with any loss function $l \in \mathcal{L}_1$ is less than $u/\sqrt{n} + \sup_{l \in \mathcal{L}_1} \mathrm{reg}_l^\theta(d_0^\theta)$ for large $n$.

PROOF OF THEOREM 5.1. By Proposition 8.1 we have

$$\sqrt{n} \sup_{l \in \mathcal{L}_1} |l^n(d_l^n) - l(\theta, d_l^\theta) - D_{10}l(\theta, d_l^\theta)^t(\theta_n - \theta)| \overset{Q^*}{\to} 0$$

and

$$\sqrt{n} \sup_{l \in \mathcal{L}_1} |l^n(d_0^n) - l(\theta, d_0^\theta)$$
$$- D_{01}l(\theta, d_0^\theta)^t(d_0^n - d_0^\theta) - D_{10}l(\theta, d_0^\theta)^t(\theta_n - \theta)| \overset{Q^*}{\to} 0.$$

The conclusion easily follows from Theorem 4.1 and Slutsky's lemma. □

From a practical point of view, it is of interest to consider the particular case where the optimal decision $d_l^\theta$ is actually independent of $l$, as is the case in estimation problems. If we assume moreover that $l_0$ is such that $d_0^\theta = d_l^\theta$, then by Theorem 5.1,

$$\sqrt{n} \sup_{l \in \mathcal{L}_1} \mathrm{reg}_l^n(d_0^n) \rightsquigarrow 0.$$

In this situation, we can expect to obtain a better rate of convergence. As a matter of fact, it turns out that the rate of convergence of the posterior regret is of order $n$.

THEOREM 5.2. *Assume that $d_0^\theta = d_l^\theta$ for every $l \in \mathcal{L}_1$. Then under the assumptions* M,

$$n \sup_{l \in \mathcal{L}_1} \mathrm{reg}_l^n(d_0^n) \rightsquigarrow \tfrac{1}{2} \sup_{l \in \mathcal{L}_1} [Z_\theta^t(\varphi(l_0) - \varphi(l))^t D_{02}l(\theta, d_l^\theta)(\varphi(l_0) - \varphi(l))Z_\theta].$$

The theorem gains in interest if we consider the special case where $\mathcal{D} = \Theta$, and $l_0$ and every $l \in \mathcal{L}_1$ are functions of $d - \sigma$, which is a very common situation in estimation problems. In this case, $\varphi(l) = -I_p$, where $I_p$ is the $p \times p$ identity matrix and

$$n \sup_{l \in \mathcal{L}_1} \mathrm{reg}_l^n(d_0^n) \rightsquigarrow 0.$$

It is easy to check that every assumption of this section is satisfied by the class of Example 4.2. Thus, the result above also holds for this class.

EXAMPLE 4.1 (continued). The assumptions of Section 5 are fulfilled with $l_0(\sigma, d) = (d - \sigma)^2$. Define $p(w, n)$ such that $\tilde{l}^n(d) = p(w, n) l^n(d)$ and assume that $\sup_{w \in \Gamma} p(w, n)$ remains bounded in $Q^*$ probability [this holds, e.g., if there exists $\underline{w}$ such that $w \geq \underline{w}$ for all $w \in \Gamma$ and if $\underline{w}$ and $w_0$ satisfy the conditions of Strasser (1976) or Abraham and Cadre (2002)]. We deduce from the above remark that $n \sup_{w \in \Gamma} (\int (d_0^n - \sigma)^2 w_n(d\sigma) - V(w_n)) \rightsquigarrow 0$, where $d_0^n$ and $V(w_n)$ are, respectively, the posterior expectation derived from the prior $w_0$ and the posterior variance derived from the prior $w$.

PROOF OF THEOREM 5.2. Since $D_{01} l^n(d_l^n) = 0$, we have, by Taylor's formula,

$$\text{reg}_l^n(d_0^n) = l^n(d_0^n) - l^n(d_l^n)$$
$$= \int_0^1 (1-s)(d_0^n - d_l^n)^t D_{02} l^n(d_l^n - s(d_0^n - d_l^n))(d_0^n - d_l^n) \, ds.$$

However, by Theorem 4.1 and Lemma 8.4,

$$\sup_{l \in \mathcal{L}_1} \sup_{s \in [0,1]} \|D_{02} l^n(d_l^n - s(d_0^n - d_l^n)) - D_{02} l(\theta, d_l^\theta)\| \overset{Q^*}{\to} 0.$$

Moreover, we easily get by Theorem 4.1 that

$$\sqrt{n}(d_0^n - d_l^n)_{l \in \mathcal{L}_1} \rightsquigarrow ((\varphi(l_0) - \varphi(l))Z_\theta)_{l \in \mathcal{L}_1}.$$

Hence

$$n \sup_{l \in \mathcal{L}_1} |\text{reg}_l^n(d_0^n) - \tfrac{1}{2}(d_0^n - d_l^n)^t D_{02} l(\theta, d_l^\theta)(d_0^n - d_l^n)| \overset{Q^*}{\to} 0.$$

We conclude by using again the asymptotic behavior of $\sqrt{n}(d_0^n - d_l^n)_{l \in \mathcal{L}_1}$. □

**6. Range of the posterior expected loss.** The beginning of this section is devoted to the study of the range of the posterior expected loss,

(6.1) $$\text{ran}_{\mathcal{S}_0}^n(d) = \sup_{l \in \mathcal{S}_0} l^n(d) - \inf_{l \in \mathcal{S}_0} l^n(d),$$

where $d \in \mathcal{D}$ and $\mathcal{S}_0$ is defined in Section 5.

THEOREM 6.1. *Assume that $d_0^\theta = d_l^\theta$ and $l(\theta, d_0^\theta) = l'(\theta, d_0^\theta)$ for every $l$ and $l' \in \mathcal{S}_0$. Then, under the assumptions* M,

$$\sqrt{n} \, \text{ran}_{\mathcal{S}_0}^n(d_0^n) \rightsquigarrow \sup_{l \in \mathcal{S}_0}[D_{10} l(\theta, d_0^\theta)^t Z_\theta] - \inf_{l \in \mathcal{S}_0}[D_{10} l(\theta, d_0^\theta)^t Z_\theta].$$

PROOF. Since $D_{01} l(\theta, d_0^\theta) = 0$, Proposition 8.1 shows that

$$\sqrt{n} \sup_{l \in \mathcal{S}_0} |l^n(d_0^n) - l(\theta, d_0^\theta) - D_{10} l(\theta, d_0^\theta)^t(\theta_n - \theta)| \overset{Q^*}{\to} 0.$$

This gives $(\sqrt{n}(l^n(d_0^n) - l(\theta, d_0^\theta)))_{l \in \mathcal{S}_0} \rightsquigarrow (D_{10}l(\theta, d_0^\theta)^t Z_\theta)_{l \in \mathcal{S}_0}$ according to Theorem 4.1, but, by assumption,

$$\mathrm{ran}_{\mathcal{S}_0}^n(d_0^n) = \sup_{l \in \mathcal{S}_0}[l^n(d_0^n) - l(\theta, d_0^\theta)] - \inf_{l \in \mathcal{S}_0}[l^n(d_0^n) - l(\theta, d_0^\theta)],$$

so that the conclusion follows from a continuity argument as in the proof of Theorem 4.1(iii).  □

It is worth pointing out that $\mathrm{ran}_{\mathcal{S}_0}^n(d) = S^n(d) - I^n(d)$ when there exist $I$ and $S$ in $\mathcal{S}_0$ such that $\sup_{l \in \mathcal{S}_0} l = S$ and $\inf_{l \in \mathcal{S}_0} l = I$. Because of the above remark, let us define another class of loss functions which is well adapted to the study of the range of posterior expected loss. Let $I \in \mathcal{S}_0$ and $S \in \mathcal{S}_0$, and define $[I, S]$ to be the class of loss functions $l \colon \Theta \times \mathcal{D} \to \mathbb{R}^+$ such that $I \leq l \leq S$. Such a class was considered in Abraham (2001). The important point to note here is that regularity assumptions are only required on $I$, $S$ and $l_0$. Thus, this class includes very irregular losses as soon as they are bounded by $I$ and $S$. This is very attractive from a practical point of view since $l_0$ can be regarded as a tractable approximation of the true loss, the accuracy of which is now given by $I$ and $S$. It is also of computational interest because it involves only two loss functions. For simplicity of notation, we write $\mathrm{ran}_{IS}^n(d)$ instead of $\mathrm{ran}_{[I,S]}^n(d)$, where the previous expression is defined by replacing $\mathcal{S}_0$ by $[I, S]$ in (6.1). Similarly, we write

$$\mathrm{ran}_{IS}^\theta(d) = \sup_{l \in [I,S]} l(\theta, d) - \inf_{l \in [I,S]} l(\theta, d).$$

THEOREM 6.2. *Under the assumptions* M,

$$\sqrt{n}(\mathrm{ran}_{IS}^n(d_0^n) - \mathrm{ran}_{IS}^\theta(d_0^\theta))$$
$$\rightsquigarrow [[D_{10}(S - I)(\theta, d_0^\theta)]^t - [D_{01}(S - I)(\theta, d_0^\theta)]^t \varphi(l_0)]Z_\theta.$$

PROOF. Since $S \in \mathcal{S}_0$, Proposition 8.1 yields

$$\sqrt{n}|S^n(d_0^n) - S(\theta, d_0^\theta) - D_{01}S(\theta, d_0^\theta)^t(d_0^n - d_0^\theta) - D_{10}S(\theta, d_0^\theta)^t(\theta_n - \theta)| \xrightarrow{Q_\theta^*} 0.$$

The same result holds with $S$ replaced by $I$. Theorem 6.2 is then an immediate consequence of Theorem 4.1 and assumption M2, since by assumption

$$\mathrm{ran}_{IS}^n(d_0^n) - \mathrm{ran}_{IS}^\theta(d_0^\theta) = [S^n(d_0^n) - S(\theta, d_0^\theta)] + [I(\theta, d_0^\theta) - I^n(d_0^n)]. \quad \square$$

Observe that if $S$, $I$ and $l_0$ are functions of $d - \sigma$, Theorem 6.2 reduces to

$$\sqrt{n}(\mathrm{ran}_{IS}^n(d_0^n) - \mathrm{ran}_{IS}^\theta(d_0^\theta)) \rightsquigarrow 0.$$

In this case we can improve the rate of convergence.

THEOREM 6.3. *Assume that $I(\cdot, d_0^\theta)$ and $S(\cdot, d_0^\theta)$ are twice continuously differentiable, $D_{10}I(\theta, d_0^\theta) = D_{10}S(\theta, d_0^\theta)$ and $D_{01}I(\theta, d_0^\theta) = D_{01}S(\theta, d_0^\theta)$. Then under the assumptions* M,

$$n(\mathrm{ran}_{IS}^n(d_0^n) - \mathrm{ran}_{IS}^\theta(d_0^\theta)) \rightsquigarrow \tfrac{1}{2}[Z_\theta^t(N_S - N_I)Z_\theta + L_S - L_I],$$

*where*

$$N_S = \varphi(l_0)^t D_{02}S(\theta, d_0^\theta)\varphi(l_0) + D_{20}S(\theta, d_0^\theta) - 2D_{11}S(\theta, d_0^\theta)^t \varphi(l_0),$$

$N_I$ *is defined by replacing $S$ by $I$ in the above formula, and the constants $L_S$ and $L_I$ are defined in Section* 8 *by replacing $f$ by $S(\cdot, d_0^\theta)$ and $I(\cdot, d_0^\theta)$, respectively, in* (8.6). *Furthermore, if $D_{01}S(\theta, d_0^\theta) = D_{10}S(\theta, d_0^\theta) = 0$, then*

(6.2) $$n(S^n(d_0^n) - S(\theta, d_0^\theta)) \rightsquigarrow \tfrac{1}{2}[Z_\theta^t N_S Z_\theta + L_S].$$

*The same result holds if $S$ is replaced by $I$ in* (6.2) *under the assumptions $D_{01}I(\theta, d_0^\theta) = D_{10}I(\theta, d_0^\theta) = 0$.*

Consider again the usual case where $l_0$, $S$ and $I$ may be expressed as functions of $d - \sigma$. Then we have $\varphi(l_0) = -I_p$, $D_{02}S = D_{20}S = -D_{11}S$ and finally $N_S = N_I = 0$, so that, by Theorem 6.3,

$$n(\mathrm{ran}_{IS}^n(d_0^n) - \mathrm{ran}_{IS}^\theta(d_0^\theta)) \rightsquigarrow \tfrac{1}{2}(L_S - L_I).$$

EXAMPLE 4.1 (continued). Take $w_I$ and $w_S$ in $\Gamma$ and consider the density ratio class $\Gamma' = \{w \in L_1(m) : w_I \leq w \leq w_S\}$. If $p(w_S, n) \overset{Q^*}{\to} w_0(\theta)/w_S(\theta)$ [which holds under the conditions of Strasser (1976) or Abraham and Cadre (2002)], it can be proved from (6.2) that

$$n \sup_{w \in \Gamma'} \int_\Theta (d_0^n - \sigma)^2 w_n(d\sigma) \left( I_\theta \int_\Theta \tau^2 F_\theta(d\tau) \right)^{-1}$$

remains asymptotically in the interval $[1, w_S(\theta)/w_I(\theta)]$.

PROOF OF THEOREM 6.3. Write $\Delta = S - I$. Let us first examine the convergence of the sequence $n(\Delta^n(d_0^n) - \Delta(\theta, d_0^\theta))$. By Taylor's formula,

$$\begin{aligned}
\Delta^n(d_0^n) &- \Delta^n(d_0^\theta) \\
&= D_{01}\Delta^n(d_0^\theta)^t(d_0^n - d_0^\theta) \\
&\quad + \int_0^1 (1-s)(d_0^n - d_0^\theta)^t D_{02}\Delta^n(d_0^\theta + s(d_0^n - d_0^\theta))(d_0^n - d_0^\theta)\, ds \\
&= A + B,
\end{aligned}$$

where $A$ and $B$ are obviously defined. Theorems 4.1 and 8.1 show that

$$n|A - (\theta_n - \theta)^t D_{11}\Delta(\theta, d_0^\theta)^t(d_0^n - d_0^\theta)| \overset{Q}{\to} 0.$$

Moreover, by Lemma 8.4 and Theorem 4.1, we have

$$n|B - \tfrac{1}{2}(d_0^n - d_0^\theta)^t D_{02}\Delta(\theta, d_0^\theta)(d_0^n - d_0^\theta)| \xrightarrow{Q} 0.$$

Finally, since $D_{10}\Delta(\theta, d_0^\theta) = 0$, Theorem 8.2 shows that

$$n|\Delta^n(d_0^\theta) - \Delta(\theta, d_0^\theta) - \tfrac{1}{2}(\theta_n - \theta)^t D_{20}\Delta(\theta, d_0^\theta)(\theta_n - \theta) - \tfrac{1}{2}L_\Delta| \xrightarrow{Q} 0.$$

Therefore, it follows from Theorem 4.1 that

$$n|\Delta^n(d_0^n) - \Delta(\theta, d_0^\theta) - \tfrac{1}{2}[(\theta_n - \theta)^t N_\Delta(\theta_n - \theta) + L_\Delta]| \xrightarrow{Q} 0.$$

The second part of Theorem 6.3 is obtained by replacing $\Delta$ by $S$ and $I$, respectively, in the above calculations.   □

## 7. Discussion.

We give in this article sufficient conditions to get optimal rates of convergence. Let us investigate whether they are necessary. We mainly discuss the existence of the second $d$ derivative.

Consider the class $\mathcal{F}$ of Section 2.1 and define a new class $\tilde{\mathcal{F}}$ by replacing $U$ and $L$, respectively, by $\tilde{U}(\sigma, d) = f(d - \sigma)$ and $\tilde{L}(\sigma, d) = f(\sigma - d)$ in the construction of $\mathcal{F}$, where $f(t) = e^{-t} + t - 1$. Note that the quadratic loss $l_0$ defined in Section 2.1 belongs to $\tilde{\mathcal{F}}$. From the arguments of Section 2.1, the diameter of $\{d_l^n, l \in \tilde{\mathcal{F}}\}$ is equal to the diameter of $\{d_{\tilde{U}}^n, d_{\tilde{L}}^n\}$. Thus, from Section 4 (Example 4.2 applied to $\mathcal{L} = \{\tilde{U}, \tilde{L}\}$), $\sqrt{n}$ diameter$\{d_l^n, l \in \tilde{\mathcal{F}}\} \rightsquigarrow 0$ while $\sqrt{n}$ diameter$\{d_l^n, l \in \mathcal{F}\} \rightsquigarrow r_1 - r_2 > 0$. The difference in the limit indicates different rates of convergence, which are due to the fact that $D_{02}U(\theta, \theta)$ does not exist while $D_{02}\tilde{U}(\sigma, \theta) \approx D_{02}\tilde{U}(\theta, \theta)$ for $\sigma$ close to $\theta$. From a technical point of view the term $D_{02}l^n(t_{l,s}^n)$, defined in the proof of Theorem 4.1, no longer converges to $D_{02}l(\theta, d_0^\theta)$ when $l = U$, but switches from $k_1$ and $k_2$ according to the sign of $d_U^n - \theta$ even for large $n$. Consequently, it is no longer possible to derive in this way the limit of $\sqrt{n}(d_U^n - \theta)$ and Theorem 4.1 does not hold for $\mathcal{L} = \{L, U\}$. [A theoretical asymptotic study of such classes can be found in Abraham (2002).] The default of smoothness [i.e., $D_{02}U(\theta, \theta)$ does not exist] slows down the rate of convergence. Analogous situations have already been noted in prior robustness: classes with point mass priors have slower rates of convergence [Sivaganesan (1988)].

## 8. Auxiliary assumptions and results.

### 8.1. *The assumptions* M.

M1. There exist $\theta \in \Theta$ and a matrix $I_\theta$ such that $\sqrt{n}(\theta_n - \theta)$ converges in distribution to a centered normal random variable $Z_\theta$ with covariance matrix $I_\theta$.

M2. For every $g \in L_1(\pi)$ and $\alpha > 0$, there exists $\eta > 0$ such that

$$e^{\eta n} \int_{\|\sigma - \theta\| \geq \alpha} g(\sigma)\pi_n(d\sigma) \to 0 \qquad \text{in } Q \text{ probability.}$$

Write for all $k > 0$,

$$W_n^k = \{\sigma \in \Theta : \|T(\sigma)\| \leq \sqrt{k \log n}\,\},$$

where $T(\sigma) = \sqrt{n} I_\theta^{-1/2}(\sigma - \theta_n)$. Let $F_n$ be the probability distribution induced by $T$ applied to $\pi_n$ and let $B_n^k$ be the closed ball with center $\theta$ and radius $\sqrt{k \log n}$.

M3. For all $r > 0$, there exist $k > 0$ and $c > 0$ such that

$$Q(\pi_n(\Theta \setminus W_n^k) > cn^{-r}) \to 0.$$

M4. There exist a probability distribution with zero mean $F_\theta$ such that

$$\int_{B_n^k} g(\sigma)F_n(d\sigma) \to \int_\Theta g(\sigma)F_\theta(d\sigma)$$

in $Q$ probability, for all $g : \Theta \to \mathbb{R}$ with $|g(\sigma)| \leq c(1 + \|\sigma\|^2)$ for some $c > 0$ and all $\sigma \in \Theta$.

### 8.2. *Asymptotics for the posterior expectation.*

Throughout this section, we denote by $Gf(\sigma)$ the gradient at $\sigma \in \Theta$ of a function $f : \Theta \to \mathbb{R}$.

#### 8.2.1. *First order result.*

We denote by $\mathcal{P}_\theta$ a set of functions $f : \Theta \to \mathbb{R}$ with the following properties:

A1. For all $f \in \mathcal{P}_\theta$, $f(\theta) = 0$.
A2. There exists an open neighborhood $V_\theta'$ of $\theta$ on which any $f \in \mathcal{P}_\theta$ is continuously differentiable and $\sup_{f \in \mathcal{P}_\theta} \|Gf(\theta)\| < \infty$.
A3. The family $\{Gf|_{V_\theta'}, f \in \mathcal{P}_\theta\}$ is equicontinuous at $\theta$.
A4. There exist a $\pi$-integrable function $q : \Theta \to \mathbb{R}$ and $\delta_0 > 0$ such that

$$\sup_{f \in \mathcal{P}_\theta} |f(\sigma)| \leq q(\sigma) \qquad \forall\, \sigma \in \Theta \quad \text{and} \quad \sup_{\|\sigma - \theta\| \leq \delta_0} q(\sigma) < \infty.$$

**THEOREM 8.1.** *Under the assumptions* M,

$$\sqrt{n} \sup_{f \in \mathcal{P}_\theta} \left| \int_\Theta f(\sigma)\pi_n(d\sigma) - Gf(\theta)^t(\theta_n - \theta) \right| \overset{Q^*}{\to} 0.$$

PROOF.  We proceed analogously to the proof of Theorem 1 of Strasser (1975). We separate the proof into two steps.

STEP 1. Let us prove that for every $c > 0$, there exists $k > 0$ such that

$$Q^*\left(\sqrt{n} \sup_{f \in \mathcal{P}_\theta} \int_{\Theta \setminus W_n^k} |f(\sigma)| \pi_n(d\sigma) > c\right) \to 0.$$

Let $i = \inf_{\|\sigma\|=1} \|I_\theta^{-1/2}\sigma\|$ and $\delta = i\delta_0$, where $\delta_0$ is the real number of A4. Clearly, we have $i > 0$ and hence $\delta > 0$. Moreover, we also have, by A4,

$$\alpha := \sup_{\|I_\theta^{-1/2}(\sigma-\theta)\| \le \delta} q(\sigma) \le \sup_{\|\sigma-\theta\| \le \delta_0} q(\sigma) < \infty.$$

Fix $c > 0$. By A4, we have, for all $k > 0$,

$$\sqrt{n} \sup_{f \in \mathcal{P}_\theta} \int_{\Theta \setminus W_n^k} |f(\sigma)| \pi_n(d\sigma) > c \quad \Longrightarrow \quad \sqrt{n} \int_{\Theta \setminus W_n^k} q(\sigma) \pi_n(d\sigma) > c,$$

and if the latter property holds, then

(8.1)
$$\|I_\theta^{-1/2}(\theta_n - \theta)\| \ge \delta/2 \qquad \text{or}$$
$$\left(\sqrt{n} \int_{\Theta \setminus W_n^k} q(\sigma) \pi_n(d\sigma) > c, \|I_\theta^{-1/2}(\theta_n - \theta)\| < \frac{\delta}{2}\right).$$

The probability of the event associated with the first property tends to 0 by M1. We now focus on the second property. Let us denote by $\mathcal{E}$ the subset of $\Theta$ defined as

$$\mathcal{E} = \{\sigma \in \Theta : \|I_\theta^{-1/2}(\sigma - \theta)\| < \delta\}.$$

There exists $N \ge 1$ such that if $\|I_\theta^{-1/2}(\theta_n - \theta)\| < \delta/2$, then for all $n \ge N$, $W_n^k \subset \mathcal{E}$. Thus, if the second property in (8.1) holds,

$$\left(\sqrt{n} \int_{\Theta \setminus \mathcal{E}} q(\sigma) \pi_n(d\sigma) > \frac{c}{2}\right) \quad \text{or} \quad \left(\sqrt{n} \int_{\mathcal{E} \setminus W_n^k} q(\sigma) \pi_n(d\sigma) > \frac{c}{2}, W_n^k \subset \mathcal{E}\right).$$

Using the obvious notation, let $A$ and $B$ be the events associated with the above properties. On one hand, the probability of $A$ tends to 0 by M2. On the other hand,

$$B \subset \{\alpha\sqrt{n}\pi_n(\Theta \setminus W_n^k) > c/2\}$$

and, for some $k > 0$, the probability of the latter event tends to 0 by M3.

STEP 2. Let us prove that for all $k, c > 0$,

$$Q^*\left(\sqrt{n} \sup_{f \in \mathcal{P}_\theta} \left|\int_{W_n^k} f(\sigma)\pi_n(d\sigma) - Gf(\theta)^t(\theta_n - \theta)\right| > c\right) \to 0.$$

We obviously have, for all $f \in \mathcal{P}_\theta$,

(8.2)
$$\int_{W_n^k} f(\sigma)\pi_n(d\sigma) = \int_{B_n^k} f(T^{-1}(\tau))F_n(d\tau),$$

where $T$ and $B_n^k$ are defined in Section 3 [recall that $T^{-1}(\tau) = \theta_n + n^{-1/2} I_\theta^{1/2} \tau$].
If $T^{-1}(\tau) \in V_\theta'$, then there exists $\lambda \in\, ]0, 1[$ such that, according to A1,

$$(8.3) \qquad f(T^{-1}(\tau)) = Gf(\theta + \lambda u(\tau))^t u(\tau),$$

where $u(\tau) = \theta_n - \theta + n^{-1/2} I_\theta^{1/2} \tau$. Let us denote by $H$ the property

$$\forall \tau \in B_n^k \qquad T^{-1}(\tau) \in V_\theta' \quad \text{and} \quad \theta + \lambda u(\tau) \in V_\theta'.$$

It is easy to check that there exist $s > 0$ and $N \geq 1$ such that, for all $n \geq N$,
$\|\theta_n - \theta\| \leq s \Longrightarrow H$. Then, if the property

$$\sqrt{n} \sup_{f \in \mathcal{P}_\theta} \left| \int_{W_n^k} f(\sigma) \pi_n(d\sigma) - Gf(\theta)^t(\theta_n - \theta) \right| > c$$

holds, we have $\|\theta_n - \theta\| > s$ or

$$(8.4) \qquad \left( \sqrt{n} \sup_{f \in \mathcal{P}_\theta} \left| \int_{W_n^k} f(\sigma) \pi_n(d\sigma) - Gf(\theta)^T(\theta_n - \theta) \right| > c, H \right).$$

By M1 we need only to focus on the latter property. If $H$ holds, we have,
according to (8.2) and (8.3),

$$\sup_{f \in \mathcal{P}_\theta} \left| \int_{W_n^k} f(\sigma) \pi_n(d\sigma) - Gf(\theta)^t(\theta_n - \theta) \right|$$

$$= \sup_{f \in \mathcal{P}_\theta} \left| \int_{B_n^k} Gf(\theta + \lambda u(\tau))^t u(\tau) F_n(d\tau) - Gf(\theta)^t(\theta_n - \theta) \right|$$

$$(8.5) \qquad \leq \sup_{f \in \mathcal{P}_\theta} \left| \int_{B_n^k} (Gf(\theta + \lambda u(\tau)) - Gf(\theta))^t u(\tau) F_n(d\tau) \right|$$

$$+ \sup_{f \in \mathcal{P}_\theta} |Gf(\theta)^t(\theta_n - \theta)(F_n(B_n^k) - 1)|$$

$$+ \sup_{f \in \mathcal{P}_\theta} n^{-1/2} \left| Gf(\theta)^t I_\theta^{1/2} \int_{B_n^k} \tau F_n(d\tau) \right|.$$

Let $\gamma > 0$. By A3 there exists $\beta > 0$ such that, for all $\sigma \in V_\theta'$ with $\|\sigma - \theta\| \leq \beta$,

$$\sup_{f \in \mathcal{P}_\theta} \|Gf(\sigma) - Gf(\theta)\| \leq \gamma.$$

Let $\bar{\alpha} = \sup_{f \in \mathcal{P}_\theta} \|Gf(\theta)\|$, which is finite by A2. For all $n \geq N$, if the property
in (8.4) holds, we have, by (8.5),

$$(\|\theta_n - \theta\| + n^{-1/2} \|I_\theta^{1/2}\| \sqrt{k \log n} > \beta),$$

$$\left( \gamma \sqrt{n} \|\theta_n - \theta\| + \gamma \|I_\theta^{1/2}\| \int_{B_n^k} \|\tau\| F_n(d\tau) > \frac{c}{3} \right),$$

$$\left( \bar{\alpha} \sqrt{n} \|\theta_n - \theta\| \ |F_n(B_n^k) - 1| > \frac{c}{3} \right)$$

or

$$\left( \bar{\alpha} \| I_\theta^{1/2} \| \left\| \int_{B_n^k} \tau F_n(d\tau) \right\| > \frac{c}{3} \right).$$

Since $F_\theta$ is centered, $\int_{B_n^k} \tau F_n(d\tau) \to 0$ in probability by M4. Hence the probability of the event associated with the latter property vanishes. The probabilities of the events related with the other properties tend to 0 by M2 and M4, for some choice $\gamma$. Step 2 is then proved and the theorem is a straightforward consequence of Steps 1 and 2.   □

8.2.2. *Second order result.* Throughout this section we denote by $Hf(\sigma)$ the Hessian matrix at $\sigma \in \Theta$ of a function $f \colon \Theta \to \mathbb{R}$ that satisfies the following properties:

B1. There exists an open neighborhood $V_\theta''$ of $\Theta$ on which $f$ is twice continuously differentiable.
B2. $f(\theta) = 0$ and $Gf(\theta) = 0$.
B3. $f$ is $\pi$-integrable.

We introduce the notation

$$(8.6) \qquad L_f = \int_\Theta (I_\theta^{1/2} \tau)^t Hf(\theta)(I_\theta^{1/2} \tau) F_\theta(d\tau),$$

provided such a quantity may be defined. Note that $F_\theta$ is normal under usual models [Strasser (1976)].

THEOREM 8.2. *Under the assumptions* M,

$$n \left| \int_\Theta f(\sigma) \pi_n(d\sigma) - \tfrac{1}{2}(\theta_n - \theta)^t Hf(\theta)(\theta_n - \theta) - \tfrac{1}{2} L_f \right| \to 0$$

*in probability.*

PROOF. Following the arguments of the first step of the proof of Theorem 8.1, we can prove that for all $c > 0$ there exists $k > 0$ such that

$$Q\left( n \int_{\Theta \setminus W_n^k} |f(\sigma)| \pi_n(d\sigma) > c \right) \to 0.$$

Hence, we need only to prove that for all $k > 0$,

$$n \left| \int_{W_n^k} f(\sigma) \pi_n(d\sigma) - \tfrac{1}{2}(\theta_n - \theta)^t Hf(\theta)(\theta_n - \theta) - \tfrac{1}{2} L_f \right| \to 0$$

in probability. We use the notation of the proof of Theorem 8.1. According to B2, if $T^{-1}(\tau) \in V_\theta''$, then there exists $\lambda \in ]0, 1[$ such that

$$(8.7) \qquad f(T^{-1}(\tau)) = \tfrac{1}{2} u(\tau)^t Hf(\theta + \lambda u(\tau)) u(\tau).$$

Fix $k > 0$ and denote by $H'$ the property

$$\forall \tau \in B_n^k \qquad T^{-1}(\tau) \in V_\theta'' \quad \text{and} \quad \theta + \lambda u(\tau) \in V_\theta''.$$

For some $s > 0$ and $N \geq 1$, we have $\|\theta_n - \theta\| \leq s \Longrightarrow H'$ for all $n \geq N$. If $H'$ holds, then according to (8.2) and (8.7),

$$
\left| \int_{W_n^k} f(\sigma) \pi_n(d\sigma) - \frac{1}{2}(\theta_n - \theta)^t H f(\theta)(\theta_n - \theta) - \frac{1}{2} L_f \right|
$$

$$
\leq \frac{1}{2} \left| \int_{B_n^k} u(\tau)^t (H f(\theta + \lambda u(\tau)) - H f(\theta)) u(\tau) F_n(d\tau) \right|
$$

(8.8)
$$
+ \frac{1}{2} \|\theta_n - \theta\|^2 \|H f(\theta)\| |F_n(B_n^k) - 1|
$$

$$
+ \frac{1}{\sqrt{n}} \|H f(\theta)\| \|\theta_n - \theta\| \|I_\theta^{1/2}\| \left\| \int_{B_n^k} \tau F_n(d\tau) \right\|
$$

$$
+ \frac{1}{2n} \left| \int_{B_n^k} (I_\theta^{1/2} \tau)^t H f(\theta)(I_\theta^{1/2} \tau) F_n(d\tau) - L_f \right|.
$$

Let $\gamma > 0$. According to B1, there exists $\beta > 0$ such that if $\sigma \in V_\theta''$ with $\|\sigma - \theta\| \leq \beta$,

$$\|H f(\sigma) - H f(\theta)\| \leq \gamma.$$

Fix $c > 0$ and let

$$L^n = \int_{B_n^k} (I_\theta^{1/2} \tau)^t H f(\theta)(I_\theta^{1/2} \tau) F_n(d\tau).$$

We deduce from (8.8) that if we have

$$n \left| \int_{W_n^k} f(\sigma) \pi_n(d\sigma) - \frac{1}{2}(\theta_n - \theta)^t H f(\theta)(\theta_n - \theta) - \frac{1}{2} L_f \right| > c,$$

then for all $n \geq N$,

$$(\|\theta_n - \theta\| > s) \quad \text{or} \quad \left( \|\theta_n - \theta\| + \frac{1}{\sqrt{n}} \|I_\theta^{1/2}\| \sqrt{k \log n} > \beta \right),$$

$$\left( \frac{n\gamma}{2} \int_{B_n^k} \|u(\tau)\|^2 F_n(d\tau) > \frac{c}{4} \right),$$

(8.9)
$$\left( \frac{n}{2} \|\theta_n - \theta\|^2 \|H f(\theta)\| |F_n(B_n^k) - 1| > \frac{c}{4} \right),$$

$$\left( \sqrt{n} \|H f(\theta)\| \|\theta_n - \theta\| \|I_\theta^{1/2}\| \left\| \int_{B_n^k} \tau F_n(d\tau) \right\| > \frac{c}{4} \right) \quad \text{or}$$

$$\left( \frac{1}{2} |L^n - L_f| > \frac{c}{4} \right).$$

According to M2 and M4, the sequence $(n \int_{B_n^k} \|u(\tau)\|^2 F_n(d\tau))_n$ is stochastically bounded and hence, for some $\gamma$, we have

$$Q\left(n \frac{\gamma}{2} \int_{B_n^k} \|u(\tau)\|^2 F_n(d\tau) > \frac{c}{4}\right) \to 0.$$

Moreover, the probability of the events associated with the other properties of (8.9) obviously vanishes according to M2 and M4. □

8.3. *Technical results for the classes $\mathcal{L}$, $\mathcal{S}$ and $\mathcal{S}_0$.*

LEMMA 8.1. *Let $g$ be a $\pi$-integrable and nonnegative real-valued function such that there exists a bounded neighborhood of $\theta$ on which $g$ is bounded. Then under the assumptions M,*

$$Q\left(\int_\Theta g(\sigma)\pi_n(d\sigma) < \infty\right) \to 1.$$

PROOF. Denote by $B$ the bounded neighborhood of $\theta$. For $t \geq 1$ let $f_n(t) = Q(\int_{B^c} g(\sigma)\pi_n(d\sigma) \geq t)$. By M2 we have

$$\sup_{t \geq 1} |f_n(t)| \leq Q\left(\int_{B^c} g(\sigma)\pi_n(d\sigma) \geq 1\right) \to 0.$$

Furthermore, $\lim_{t \to \infty} f_n(t)$ exists since $f_n$ is decreasing and bounded, so that $\lim_{n \to \infty} \lim_{t \to \infty} f_n(t) = \lim_{t \to \infty} \lim_{n \to \infty} f_n(t) = 0$. We conclude by noting that

$$Q\left(\int_\Theta g(\sigma)\pi_n(d\sigma) = \infty\right) \leq \lim_{t \nearrow \infty} Q\left(\int_B g(\sigma)\pi_n(d\sigma) \geq t\right) + \lim_{t \nearrow \infty} f_n(t),$$

hence the lemma, since $g$ is bounded on $B$. □

LEMMA 8.2. *Under the assumptions M, there exists a compact set $K \subset \mathcal{D}$ such that $Q^*(\exists\, l \in \mathcal{L},\ d_l^n \in K^c) \to 0$.*

PROOF. Take $r > 0$ and $K$ compact as in 1f and introduce $\alpha$ and $0 < \varepsilon < 1$ such that

$$\sup_{l \in \mathcal{L}} \inf_{d \in K} l(\theta, d) < (1-\varepsilon)\alpha < \alpha < \inf_{l \in \mathcal{L}} \inf_{\sigma \in B(\theta, r)} \inf_{d \in K^c} l(\sigma, d).$$

Then, if $d \in K^c$, we have

$$l^n(d) = \int_{B(\theta, r)} l(\sigma, d)\pi_n(d\sigma) + \int_{B^c(\theta, r)} l(\sigma, d)\pi_n(d\sigma)$$
$$> \alpha \pi_n(B(\theta, r)).$$

Thus

$$\exists l \in \mathcal{L},\, d_l^n \in K^c$$

$$\Longrightarrow \quad \exists l \in \mathcal{L},\, \exists d \in K^c,\, l^n(d) \leq \inf_{t \in K} l^n(t)$$

$$\Longrightarrow \quad \exists l \in \mathcal{L},\, \alpha \pi_n(B(\theta, r)) < \inf_{t \in K} l^n(t)$$

$$\Longrightarrow \quad \left( \exists l \in \mathcal{L},\, \alpha \pi_n(B(\theta, r)) < \inf_{t \in K} l(\theta, t) + \varepsilon \frac{\alpha}{2} \right) \quad \text{or}$$

$$\left( \exists l \in \mathcal{L},\, \inf_{t \in K} l^n(t) > \inf_{t \in K} l(\theta, t) + \varepsilon \frac{\alpha}{2} \right)$$

$$\Longrightarrow \quad \left( \alpha \left( -\frac{\varepsilon}{2} + \pi_n(B(\theta, r)) \right) \leq \sup_{l \in \mathcal{L}} \inf_{t \in K} l(\theta, t) \right) \quad \text{or}$$

$$\left( \sup_{l \in \mathcal{L}} \sup_{t \in K} |l^n(t) - l(\theta, t)| > \varepsilon \frac{\alpha}{2} \right)$$

$$\Longrightarrow \quad \left( \pi_n(B^c(\theta, r)) \geq \frac{\varepsilon}{2} \right) \quad \text{or}$$

$$\left( \int \sup_{l \in \mathcal{L}} \sup_{t \in K} |l(\sigma, t) - l(\theta, t)| \pi_n(d\sigma) > \varepsilon \frac{\alpha}{2} \right).$$

By M2, $Q(\pi_n(B^c(\theta, r)) \geq \varepsilon/2) \to 0$. Moreover, if the last condition on the right-hand side holds, then for all $\rho > 0$,

$$\int_{B(\theta, \rho)} \sup_{l \in \mathcal{L}} \sup_{t \in K} |l(\sigma, t) - l(\theta, t)| \pi_n(d\sigma) > \varepsilon \frac{\alpha}{4}$$

or

$$\int_{B^c(\theta, \rho)} \sup_{l \in \mathcal{L}} \sup_{t \in K} |l(\sigma, t) - l(\theta, t)| \pi_n(\sigma) > \varepsilon \frac{\alpha}{4}.$$

By 1d we choose $\rho$ small enough so that the outer probability of the event associated with the first property tends to 0. Then, for the second property, bound the integrand by $g_1 \in L_1(\pi)$ and conclude by the concentration assumption M2. Since $\mathcal{L}$ is $\pi$-dominated, the existence of $g_1$ is deduced from the compactness of $K$. $\quad \square$

LEMMA 8.3. *Under the assumptions* M,

$$\sup_{l \in \mathcal{L}} \|d_l^n - d_l^\theta\| \xrightarrow{Q^*} 0.$$

PROOF.   According to Lemma 8.2, we may restrict our attention to those $x \in \{x \in \mathcal{X},\, \forall l \in \mathcal{L},\, d_l^n(x) \in K\}$, where $K$ is a compact set. By 1f there is no loss of generality in assuming that $d_l^\theta \in K$ for $l \in \mathcal{L}$. Let $\varepsilon > 0$. Note that, for $l \in \mathcal{L}$ and $d \in B^c(d_l^\theta, \varepsilon)$, the property $l^n(d) \leq l^n(d_l^\theta)$ implies that

$$
\begin{aligned}
l^n(d) - l(\theta, d) &\leq -(l(\theta, d) - l(\theta, d_l^\theta)) + (l^n(d_l^\theta) - l(\theta, d_l^\theta)) \\
&\leq -\kappa(\varepsilon) + (l^n(d_l^\theta) - l(\theta, d_l^\theta)),
\end{aligned}
$$

where the last inequality follows from 1g. According to the above remark, we have, for all $r > 0$,

$$
\sup_{l \in \mathcal{L}} \|d_l^n - d_l^\theta\| > \varepsilon
$$

$$
\implies \exists\, l \in \mathcal{L},\, \exists\, d \in B^c(d_l^\theta, \varepsilon) \cap K,\, l^n(d) \leq l^n(d_l^\theta)
$$

$$
\implies \left( \sup_{l \in \mathcal{L}} \sup_{d \in B^c(d_l^\theta, \varepsilon) \cap K} |l^n(d) - l(\theta, d)| \geq \frac{\kappa(\varepsilon)}{2} \right) \quad \text{or}
$$

$$
\left( \sup_{l \in \mathcal{L}} |l^n(d_l^\theta) - l(\theta, d_l^\theta)| \geq \frac{\kappa(\varepsilon)}{2} \right)
$$

$$
\implies \sup_{l \in \mathcal{L}} \sup_{d \in K} |l^n(d) - l(\theta, d)| \geq \frac{\kappa(\varepsilon)}{2}
$$

$$
\implies \sup_{l \in \mathcal{L}} \sup_{d \in K} \sup_{\sigma \in B(\theta, r)} |l(\sigma, d) - l(\theta, d)|
$$

$$
+ \int_{B^c(\theta, r)} \sup_{l \in \mathcal{L}} \sup_{d \in K} |l(\sigma, d) - l(\theta, d)| \pi_n(d\sigma) \geq \frac{\kappa(\varepsilon)}{2}.
$$

By 1d, we can choose $r > 0$ such that

$$
\sup_{l \in \mathcal{L}} \sup_{d \in K} \sup_{\sigma \in B(\theta, r)} |l(\sigma, d) - l(\theta, d)| < \frac{\kappa(\varepsilon)}{4}
$$

and we thus get

$$
\sup_{l \in \mathcal{L}} \|d_l^n - d_l^\theta\| > \varepsilon \quad \implies \quad \int_{B^c(\theta, r)} \sup_{l \in \mathcal{L}} \sup_{d \in K} |l(\sigma, d) - l(\theta, d)| \pi_n(d\sigma) \geq \frac{\kappa(\varepsilon)}{4}.
$$

Taking into account the compactness of $K$, we can deduce from the definition of a locally $\pi$-dominated class that there exists $g_1 \in L_1(\pi)$ such that

$$
\sup_{l \in \mathcal{L}} \sup_{d \in K} |l(\sigma, d) - l(\theta, d)| \leq g_1(\sigma) \qquad \forall\, \sigma \in \Theta.
$$

The conclusion then follows from M2 and 1g. $\quad \square$

LEMMA 8.4.  *For every $n \geq 1$, $s \in [0, 1]$ and $l \in \mathcal{L}$, let $t_{l,s}^n : \mathcal{X} \to \mathcal{D}$ be a map such that $\sup_{l \in \mathcal{L}} \sup_{s \in [0,1]} \|t_{l,s}^n - d_l^\theta\| \xrightarrow{Q^*} 0$. Then, under the assumptions* M,

$$
\sup_{l \in \mathcal{L}} \sup_{s \in [0,1]} \|D_{02} l^n(t_{l,s}^n) - D_{02} l(\theta, d_l^\theta)\| \xrightarrow{Q^*} 0.
$$

PROOF. Fix $\varepsilon > 0$. By 1e take $\eta > 0$ such that $\sup_{\sigma \in V_\theta} \rho_\eta(\sigma) < \varepsilon/2$. Then

$$\sup_{l \in \mathcal{L}} \sup_{s \in [0,1]} \|D_{02} l^n(t_{l,s}^n) - D_{02} l^n(d_l^\theta)\| > \varepsilon$$

$$\implies \quad \left( \sup_{l \in \mathcal{L}} \sup_{s \in [0,1]} \|t_{l,s}^n - d_l^\theta\| > \eta \right) \quad \text{or}$$

$$\left( \int_{V_{\theta^c}} \rho_\eta(\sigma) \pi_n(d\sigma) > \frac{\varepsilon}{2} \right).$$

The outer probability of the events associated with the above properties tends to 0 by assumption and M2. Consequently, it remains to prove that

$$\sup_{l \in \mathcal{L}} \sup_{s \in [0,1]} \|D_{02} l^n(d_l^\theta) - D_{02} l(\theta, d_l^\theta)\| \xrightarrow{Q_*^*} 0.$$

By 1d take $\beta > 0$ such that

$$\sup_{l \in \mathcal{L}} \sup_{\sigma \in B(\theta, \beta)} \|D_{02} l(\sigma, d_l^\theta) - D_{02} l(\theta, d_l^\theta)\| \leq \frac{\varepsilon}{2}.$$

Then by splitting the integral according to $\Theta = B(\theta, \beta) \cup B(\theta, \beta)^c$, we have

$$\sup_{l \in \mathcal{L}} \left\| \int_\Theta (D_{02} l(\sigma, d_l^\theta) - D_{02} l(\theta, d_l^\theta)) \pi_n(d\sigma) \right\| > \varepsilon$$

$$\implies \int_{B(\theta, \beta)^c} \sup_{l \in \mathcal{L}} (\|D_{02} l(\sigma, d_l^\theta)\| + \|D_{02} l(\theta, d_l^\theta)\|) \pi_n(d\sigma) > \frac{\varepsilon}{2}.$$

Taking into account that $\mathcal{L}$ is locally $\pi$-dominated, the outer probability of the event associated with the above property tends to 0 by M2. $\square$

Following the arguments of the proof of Lemma 8.4, we obtain the result below.

LEMMA 8.5. *For every $n \geq 1$, $s \in [0,1]$ and $l \in \mathcal{S}_0$, let $t_{l,s}^n : \mathcal{X} \to \mathcal{D}$ be a map such that $\sup_{l \in \mathcal{S}_0} \sup_{s \in [0,1]} \|t_{l,s}^n - d_0^\theta\| \xrightarrow{Q_*^*} 0$. Then, under the assumptions* M,

$$\sup_{l \in \mathcal{S}_0} \sup_{s \in [0,1]} \|D_{01} l^n(t_{l,s}^n) - D_{01} l(\theta, d_0^\theta)\| \xrightarrow{Q_*^*} 0.$$

*The result is still true if $d_0^\theta$ and $\mathcal{S}_0$ are replaced by $d_l^\theta$ and $\mathcal{S}$, respectively, in which case $D_{01} l(\theta, d_l^\theta) = 0$.*

PROPOSITION 8.1. *Under the assumptions* M,

$$\sqrt{n} \sup_{l \in \mathcal{S}_0} |l^n(d_0^n) - l(\theta, d_0^\theta)$$

$$- D_{01} l(\theta, d_0^\theta)^t (d_0^n - d_0^\theta) - D_{10} l(\theta, d_0^\theta)^t (\theta_n - \theta)| \xrightarrow{Q_*^*} 0.$$

*The result is still true if $d_0^\theta$, $d_0^n$ and $\mathcal{S}_0$ are replaced by $d_l^\theta$, $d_l^n$ and $\mathcal{S}$, respectively, in which case $D_{01}l(\theta, d_l^\theta) = 0$.*

Proof. Let $l \in \mathcal{S}_0$. Then

$$l^n(d_0^n) - l(\theta, d_0^\theta) = (l^n(d_0^n) - l^n(d_0^\theta)) + (l^n(d_0^\theta) - l(\theta, d_0^\theta)).$$

By Taylor's formula, the first term on the right-hand side equals

$$\int_0^1 D_{01}l^n(d_0^\theta + s(d_0^n - d_0^\theta))^t (d_0^n - d_0^\theta) \, ds,$$

so that, by Lemma 8.5 and Theorem 4.1,

$$\sqrt{n} \sup_{l \in \mathcal{S}_0} |l^n(d_0^n) - l^n(d_0^\theta) - D_{01}l(\theta, d_0^\theta)^t (d_0^n - d_0^\theta)| \overset{Q^*}{\to} 0.$$

Moreover, by Theorem 8.1,

$$\sqrt{n} \sup_{l \in \mathcal{S}_0} |l^n(d_0^\theta) - l(\theta, d_0^\theta) - D_{10}l(\theta, d_0^\theta)^t (\theta_n - \theta)| \overset{Q^*}{\to} 0,$$

which proves the proposition. $\square$

8.4. *Technical result related to weak convergence.* Let $F(\mathcal{L})$ be the set of mappings from $\mathcal{L}$ into $\mathbb{R}$ and let $\mathcal{M}_{i,j}(F(\mathcal{L}))$ be the set of $i \times j$ matrices with coefficient in $F(\mathcal{L})$. For $A \in \mathcal{M}_{i,j}(F(\mathcal{L}))$, write $\|A\|_\infty = \sup_{l \in \mathcal{L}} \|A(t)\|$. The proof of the following lemma is left to the reader.

Lemma 8.6. *For all $n \geq 1$, consider the maps $M_n : \mathcal{X} \to \mathcal{M}_{p,1}(F(\mathcal{L}))$, $A_n : \mathcal{X} \to \mathcal{M}_{p,p}(F(\mathcal{L}))$ and $R_n : \mathcal{X} \to \mathcal{M}_{p,1}(F(\mathcal{L}))$. Let $A \in \mathcal{M}_{p,p}(F(\mathcal{L}))$ such that $\inf_{l \in \mathcal{L}} |\det A(l)| > 0$ and $\|A\| < \infty$. Assume that $A_n \overset{Q^*}{\to} A$, $R_n \rightsquigarrow R$, where $R : \mathcal{X} \to \mathcal{M}_{p,1}(F(\mathcal{L}))$ is Borel measurable and $\|A_n M_n - R_n\|_\infty \overset{Q^*}{\to} 0$. Then we have $\|M_n - A^{-1}R_n\|_\infty \overset{Q^*}{\to} 0$.*

**Acknowledgments.** We thank the referees and the Associate Editor for their helpful comments.

ENSAM–INRA
UMR Biométrie et Analyse de Systèmes
2 Place P. Viala
34060 Montpellier Cedex 1
France
e-mail: abraham@helios.ensam.inra.fr

Départment de Probabilités
et Statistiques
Université Montpellier II
CC051 Place E. Battaillon
34095 Montpellier Cedex 5
France